\documentclass[12pt]{amsart}
\usepackage{amsmath,amssymb,latexsym}
\def\dcl{\mathrm{dcl}}
\def\acl{\mathrm{acl}}
\def\cl{\mathrm{cl}}
\def\bdd{\mathrm{bdd}}
\def\cb{\mathrm{Cb}}

\def\tp{\mathrm{tp}}
\def\lstp{\mathrm{lstp}}
\def\domeq{\mathbin{\underline\square}}
\def\Ind#1#2{#1\setbox0=\hbox{$#1x$}\kern\wd0\hbox to 0pt{\hss$#1\mid$\hss}
\lower.9\ht0\hbox to 0pt{\hss$#1\smile$\hss}\kern\wd0}
\def\ind{\mathop{\mathpalette\Ind{}}}
\def\Notind#1#2{#1\setbox0=\hbox{$#1x$}\kern\wd0\hbox to 0pt{\mathchardef
\nn="3236\hss$#1\nn$\kern1.4\wd0\hss}\hbox to 0pt{\hss$#1\mid$\hss}\lower.9\ht0
\hbox to 0pt{\hss$#1\smile$\hss}\kern\wd0}
\def\nind{\mathop{\mathpalette\Notind{}}}

\theoremstyle{plain}
\newtheorem{theorem}{Theorem}[section]
\newtheorem{proposition}[theorem]{Proposition}
\newtheorem{fact}[theorem]{Fact}
\newtheorem{lemma}[theorem]{Lemma}
\newtheorem{corollary}[theorem]{Corollary}
\newtheorem*{claim}{Claim}

\theoremstyle{definition}
\newtheorem{definition}[theorem]{Definition}
\newtheorem{remark}[theorem]{Remark}
\newtheorem{expl}[theorem]{Example}

\def\bsp{\begin{expl}}
\def\ebsp{\end{expl}}
\def\beh{\begin{claim}}
\def\ebeh{\end{claim}}
\def\defn{\begin{definition}}
\def\edefn{\end{definition}}
\def\satz{\begin{theorem}}
\def\esatz{\end{theorem}}
\def\tats{\begin{fact}}
\def\etats{\end{fact}}
\def\kor{\begin{corollary}}
\def\ekor{\end{corollary}}
\def\lmm{\begin{lemma}}
\def\elmm{\end{lemma}}
\def\bem{\begin{remark}}
\def\ebem{\end{remark}}
\def\bew{\par\noindent{\em Proof: }}
\def\bewbeh{\par\noindent{\em Proof of Claim: }}
\def\satzli{\begin{proposition}}
\def\esatzli{\end{proposition}}

\begin{document}
\title{Plus Ultra}
\author{Frank O. Wagner}
\address{Universit\'e de Lyon; CNRS; Universit\'e Lyon 1; Institut Camille Jordan UMR5208, 43 bd du 11 novembre 1918, 69622 Villeurbanne Cedex, France}
\email{wagner@math.univ-lyon1.fr}
\date{21 March 2014}
\thanks{Partially supported by ANR-09-BLAN-0047 Modig and ANR-13-BS01-0006}
\begin{abstract}We define a reasonably well-behaved class of ultraimaginaries,
i.e.\ classes modulo $\emptyset$-invariant equivalence relations, called {\em
tame}, and establish some basic simplicity-theoretic facts. We also show feeble
elimination of supersimple ultraimaginaries: If $e$ is an ultraimaginary
definable over a tuple $a$ with $SU(a)<\omega^{\alpha+1}$, then $e$ is
eliminable up to rank $<\omega^\alpha$. Finally, we prove some uniform versions
of the weak canonical base property.\end{abstract}
\maketitle

\section{Introduction}

This paper arose out of an attempt to understand and generalize Chatzidakis'
results on the weak canonical base property \cite[Proposition 1.14 and Lemma
1.15]{zoe}. In doing so, we realized that certain stability-theoretic phenomena
were best explained using ultraimaginaries. It should be noted that
ultraimaginaries occur naturally in simplicity theory and were in fact briefly
considered in \cite{bytw05} before specializing to the more restricted class of
almost hyperimaginaries. However, they have faded into oblivion since Ben Yaacov
\cite{BY03} has shown that no satisfactory independence theory can exist for
them, as there are problems with both the finite character and the
extension axiom for independence. Nevertheless, at least finite character can be
salvaged if one restricts to quasi-finitary ultraimaginaries in a supersimple
theory, or more generally to what we call {\em tame} ultraimaginaries. 

We shall define ultraimaginaries in Section \ref{sec2} and give various examples. We also give a first example of a natural general result involving them, Proposition \ref{lascar}, which for a supersimple theory of finite rank specializes to a theorem of Lascar. In Section \ref{sec3} we define tame ultraimaginaries and recover certain tools from simplicity theory, even though, due to the lack of extension, canonical bases are not available in our context. One may thus hope to extend the techniques of this section for instance to the superrosy context, where the lack of canonical bases has been one of the main technical problems.

In Section \ref{sec4} we prove feeble elimination of ultraimaginaries. In
particular ultraimaginaries of finite rank are interbounded with
hyperimaginaries. This is used in Section \ref{sec5} to generalize some of
Chatzidakis' results \cite{zoe} on the weak canonical base property from sets of
finite $SU$-rank to arbitrary ordinal $SU$-rank. It is interesting to compare
this generalization to the coarser \cite[Theorem 5.4]{pw} which uses
$\alpha$-closure. We expect that this is a general phenomenon: The use of
ultraimaginaries allows for a more direct and more refined proof without
explicit use of $SU$-rank; rank considerations principally intervene via the
feeble elimination result of Section \ref{sec4} and the technical results of
Section \ref{sec3}.

All elements, tuples and parameter sets are hyperimaginary, unless stated
otherwise. For an introduction to simplicity and hyperimaginaries, the reader is
invited to consult \cite{ca11}, \cite{ki14} or \cite{wa00}.

\section{Ultraimaginaries}\label{sec2}

\defn An {\em ultraimaginary} is the class $a_E$ of a tuple $a$ under an
$\emptyset$-invariant equivalence relation $E$.\edefn
Note that tuples of ultraimaginaires are again ultraimaginary.
\defn An ultraimaginary $a_E$ is {\em definable} over an ultraimaginary $b_F$
if any automorphism of the monster model stabilising the $F$-class of $b$
also stabilises the $E$-class of $a$. It is {\em bounded} over $b_F$ if the
orbit of $a$ under the group of automorphisms of the monster model which
stabilise the $F$-class of $b$ is contained in boundedly many
$E$-classes. A {\em representative} of an ultraimaginary $e$ is any
hyperimaginary $a$ such that $e$ is definable over $a$. An ultraimaginary
is {\em finitary\/} if it has a finite
real representative. Two (tuples of) ultraimaginaries are {\em equivalent}
over some set $A$ of parameters if they are interdefinable over $A$.\edefn
Note that in contrast to hyperimaginaries, the class of a tuple of size $\kappa$ modulo an $\emptyset$-invariant equivalence relation need not be equivalent to a tuple of ultraimaginaries with representatives of smaller size: Consider the equivalence relation on sequences of length $\kappa$ of being equal except for a subsequence of smaller length.
\bem As usual, if $E_A(x,y)$ is an $A$-invariant equivalence relation, one considers the $\emptyset$-invariant relation
$E(xX,yY)$ given by
$$(X=Y\land X\equiv A\land E_X(x,y))\lor (X=Y\land x=y).$$
This is an equivalence relation, and $(aA)_E$ is equivalent to $a_{E_A}$ over
$A$.\ebem
\bem As any $\emptyset$-invariant relation, $E$ is given by a union of types over $\emptyset$.\ebem
\defn We shall say that two ultraimaginaries {\em have the same (Lascar strong)
type} over some set $A$ if they have representatives which do. If the
ambient theory is simple, we call two ultraimaginaries {\em independent} over
$A$ if they have representatives which are.\edefn
Clearly, two ultraimaginaries are conjugate by a (Lascar strong) automorphism
over $A$ if and only if they have the same (Lascar strong) type over $A$.
\bem\label{finchar} If $e$ or $e'$ is a sequence of ultraimaginaries, for
$e\ind e'$ to hold we require sequences of representatives which are
independent. In particular, it is not
clear even for real $e'$ that an infinite sequence $e$ of ultraimaginaries is
independent of $e'$ if every finite subsequence is independent of $e'$.\ebem
Ultraimaginaries arise quite naturally in stability and simplicity theory.
\bsp Let $p_A\in S(A)$ be a regular type in a stable theory. For $A',A''\models\tp(A)$ put $E(A',A'')$ if $p_{A'}\not\perp p_{A''}$. Then $E$ is an $\emptyset$-invariant equivalence relation, and $A_E$ codes the non-orthogonality class of $p_A$.\ebsp
The work with ultraimaginaries requires caution, as some basic properties become problematic.
\bsp\label{bspBY}\cite{BY03} Let $E$ be the $\emptyset$-invariant equivalence
relation on infinite sequences which holds if they differ only on finitely many
elements. Consider a sequence $I=(a_i:i<\omega)$ of elements such that no finite
subtuple is bounded over the remaining elements. Then every finite tuple $\bar
a\in I$ can be moved to a disjoint conjugate over $I_E$, but $I$ cannot.
Similarly, if $I$ is a Morley sequence in a simple theory, then $\bar a\ind I_E$
for any finite $\bar a\in I$, but $I\nind I_E$.\ebsp
Even in the $\omega$-stable context, for classes of finite tuples, the theory is not smooth.
\bsp\label{forest} Let $T$ be the theory of a cycle-free graph (forest) of infinite valency, with predicates $P_n(x,y)$ for couples of points of distance $n$ for all $<\omega$. It is easy to see by back-and-forth that $T$ eliminates quantifiers and is $\omega$-stable of rank $\omega$; the formula $P_n(x,a)$ has rank $n$ over $a$. Let $E$ be the $\emptyset$-invariant equivalence relation of being in the same connected component. Then existence of non-forking extensions fails over $a_E$, as any two points in the connected component of $a$ have some finite distance $n$, and hence rank $n$ over one another, but rank $\ge k$ over $a_E$ for all $k<\omega$, since $a_E$ is definable over any point of distance at least $k$.\par
The same phenomenon can be observed for any type $p$ of rank $SU(p)=\omega$ in a
simple theory, with the relation $E(x,y)$ on $p$ which holds if $SU(x/y)<\omega$
and $SU(y/x)<\omega$ (actually, one follows from the other by the Lascar
inequalities).\ebsp
The behaviour of Example \ref{forest} is inconvenient and signifies that we
shall avoid considering types {\em over} an ultraimaginary. The behaviour of
Example \ref{bspBY} is outright vexatious; we shall restrict the class of
ultraimaginaries under consideration in order to preserve the finite character
of independence.
\defn An ultraimaginary $e$ is {\em quasi-finitary} if there is a finite
real tuple $a$ such that $e$ is bounded over $a$.\edefn

For hyperimaginary tuples contained in the bounded closure of a finite set, we
shall use {\em quasi-finite} rather than {\em quasi-finitary}, in order to
emphasize the distinction between usual hyperimaginaries and ultraimaginaries.
The set of all (quasi-finitary) ultraimaginaries definable over
some ultraimaginary set $E$ will be denoted by $\dcl^u(E)$ (or $\dcl^{qfu}(E)$,
respectively). Similarly, $\bdd^u(E)$ and $\bdd^{qfu}(E)$ will denote the
corresponding bounded closures. If $A$
is a set of representatives for $E$, then the number of ultraimaginaries with
representatives of length $\kappa$ in $\bdd^u(E)$ (and {\em a fortiori} in
the other closures as well) is bounded in terms of the number of Lascar strong
types over $A$ of real tuples of length $\kappa$, since equality of Lascar
strong type over $A$ is the finest bounded $A$-invariant equivalence relation.
\bem\label{qfu} If $e$ is a quasi-finitary ultraimaginary then
$\bdd^u(e)=\bdd^{qfu}(e)$ and $\dcl^u(e)=\dcl^{qfu}(e)$.\ebem
\satzli\label{lascar} The following are equivalent for two
ultraimaginaries $a$ and $b$:\begin{enumerate}
\item $\bdd^u(a)\cap\bdd^u(b)=\bdd^u(\emptyset)$.
\item For any $a'\equiv^{lstp}a$ there is $n<\omega$ and a
sequence $(a_ib_i:i\le n)$ such that $$a_0=a,\quad b_0=b,\quad a_n=a'$$
and for each $i<n$
$$\bdd^u(a_i)b_{i+1}\equiv\bdd^u(a_i)b_i\quad\mbox{and}\quad
a_{i+1}\bdd^u(b_{i+1})\equiv a_i\bdd^u(b_{i+1}).$$\end{enumerate}
If $a$ or $b$ is quasi-finite, this is also eqivalent to
$\bdd^{qfu}(a)\cap\bdd^{qfu}(b)=\bdd^{qfu}(\emptyset)$.\esatzli
\bew $(1)\Rightarrow(2)$ Suppose $\bdd^u(a)\cap\bdd^u(b)=\bdd^u(\emptyset)$, and
define an $\emptyset$-invariant equivalence relation on $\lstp(ab)$ by
$E(xy,x'y')$ if there is a sequence $(x_iy_i:i\le n)$ such that $x_0y_0=xy$,
$x_ny_n=x'y'$, and for each $i<n$ we have
$$\bdd^u(x_i)y_{i+1}\equiv\bdd^u(x_i)y_i\quad\mbox{and}\quad
x_{i+1}\bdd^u(y_{i+1})\equiv x_i\bdd^u(y_{i+1}).$$
Now if $\bdd^u(a)b'\equiv\bdd^u(a)b$, then $\models
E(ab,ab')$. Hence $(ab)_E\in\bdd^u(a)$. Similarly $(ab)_E\in\bdd^u(b)$, whence
$(ab)_E\in\bdd^u(\emptyset)$. But for any $a'\equiv^{lstp}a$ there
is $b'$ with $ab\equiv^{lstp}a'b'$. Then $\models E(ab,a'b')$, in particular
$(2)$ holds.

$(2)\Rightarrow(1)$ Suppose not, and consider
$e\in(\bdd^u(a)\cap\bdd^u(b))\setminus\bdd^u(\emptyset)$. As
$e\notin\bdd^u(\emptyset)$ there is $a'\models\lstp(a)$ with
$e\notin\bdd^u(a')$. Consider a sequence $(a_i,b_i:i\le n)$ as in $(2)$. Since
$\bdd^u(a_i)b_{i+1}\equiv\bdd^u(a_i)b_i$ and
$a_{i+1}\bdd^u(b_{i+1})\equiv a_i\bdd^u(b_{i+1})$ we have
\begin{align*}\bdd^u(a_i)\cap\bdd^u(b_i)&=\bdd^u(a_i)\cap\bdd^u(b_{i+1})\\
 &=\bdd^u(a_{i+1})\cap\bdd^u(b_{i+1}).\end{align*}
In particular,
$$\begin{aligned}e\in\bdd^u(a)\cap\bdd^u(b)&=\bdd^u(a_0)\cap\bdd^u(b_0)\\
&=\bdd^u(a_n)\cap\bdd^u(b_n)\subseteq\bdd^u(a'),\end{aligned}$$
a contradiction.

The last assertion follows from Remark \ref{qfu}.\qed
\bem\label{lstp} For hyperimaginary $a$ and ultraimaginary $b$ and $c$ the
condition $\bdd^u(a)b\equiv\bdd^u(a)c$ is equivalent to
$b\equiv^{lstp}_ac$.\ebem
Using weak elimination of ultraimaginaries proven in Section \ref{sec4}, we recover a Lemma of Lascar \cite{L} (see also \cite[Lemma 2.2]{MP}), proved originally for stable theories of finite Lascar rank.
\kor Let $T$ be a simple theory of finite $SU$-rank and
$a,b$ finite imaginary tuples. The following are
equivalent:\begin{enumerate}
\item $\acl^{eq}(a)\cap\acl^{eq}(b)=\acl^{eq}(\emptyset)$.
\item For any $a'\equiv^{lstp}a$ independent of $a$
there are sequences $a=a_0,\ldots,a_n=a'$ and $b=b_0,\ldots,b_n$, such that
$b_{i+1}\equiv^{lstp}_{a_i}b_i$ and
$a_{i+1}\equiv^{lstp}_{b_{i+1}}a_i$ for each
$i<n$.\end{enumerate}\ekor
\bew By Theorem \ref{EUI} supersimple theories of
finite rank have weak elimination of quasi-finitary ultraimaginaries; by
\cite{BPW01} supersimple theories eliminate hyperimaginaries. Hence condition
(1) is equivalent to $\bdd^u(a)\cap\bdd^u(b)=\bdd^u(\emptyset)$. By Remark
\ref{lstp} condition (2) is equivalent to condition (2) of Proposition
\ref{lascar}. So $(1)\Rightarrow(2)$ follows from Proposition \ref{lascar}; for
the converse given arbitrary $a'\equiv^{lstp}_Aa$ we consider
$a''\equiv^{lstp}_Aa$ with $a''\ind_Aaa'$ and compose the sequence
$(a_ib_i:i\le n)$ from $ab=a_0b_0$ to $a_n=a''$ with the sequence $(a_ib_i:n\le
i\le\ell)$ from $a_nb_n$ to $a_\ell=a'$. Hence (2) holds for arbitrary
$a'\equiv^{lstp}_Aa$, so we can again apply Proposition
\ref{lascar}.\qed

\section{Ultraimaginaries in simple theories}\label{sec3}
From now on the ambient theory will be simple. Our notation is standard and
follows \cite{wa00}. We shall be working in  a sufficiently saturated model of
the ambient theory. Tuples are again tuples of hyperimaginaries, and closures
(definable and bounded closures) will include hyperimaginaries.
\bem\label{LKP} Since in a simple theory Lascar strong type equals Kim-Pillay strong type, we have $\bdd^u(A)=\dcl^u(\bdd(A))$. But of course, as with real and imaginary algebraic closures, $\bdd(A)\cap\bdd(B)=\bdd(\emptyset)$ does not imply $\bdd^u(A)\cap\bdd^u(B)=\bdd^u(\emptyset)$ unless the theory weakly eliminates ultraimaginaries.\ebem
In a simple theory, ultraimaginary independence is clearly symmetric, and
satisfies local character and extension (but recall that we only consider
hyperimaginary base sets), since this is inherited from suitable
representatives. As for transitivity, we have the following.
\tats\label{transitive}\cite[Lemma 1.10]{bytw05} Let $A,a$ be hyperimaginary,
and $e,e',e''$ ultraimaginary.\begin{itemize}
\item If $e\ind_A e'e''$ and $e'\ind_A e''$, then $ee'\ind_A e''$ and $e \ind_A e'$.
\item $e \ind_A ae'$ if and only if $e \ind_A a$ and $e \ind_{Aa} e'$.
\end{itemize}\etats
The Independence Theorem and Boundedness axiom also hold.
\tats\cite[page 189]{bytw05} Let $A$ be hyperimaginary and $e,e'$ ultraimaginary with $e\ind_A e'$.\begin{itemize}
\item If $f,f'$ are ultraimaginary with $f\ind_Ae$, $f'\ind_Ae'$ and $f\equiv^{Lstp}_Af'$, then there is $f''\ind_Aee'$ with $ef''\equiv_A^{Lstp}ef$ and $e'f''\equiv_A^{Lstp}e'f'$.
\item If $e''\in\bdd^u(Ae)$ then $e''\ind_A e'$. Moreover, if $e\ind_a e$ for every representative $a$ of an ultraimaginary $e''$, then $e\in\bdd^u(e'')$.\end{itemize}\etats
Next, ultraimaginary bounded closures of independent sets intersect trivially.
\lmm\label{indepintersection} If $A$ is hyperimaginary and $e,e'$ ultraimaginary with $e\ind_Ae'$, then $\bdd^u(Ae)\cap\bdd^u(Ae')=\bdd^u(A)$.\elmm
\bew Replacing $e$ and $e'$ by $A$-independent representatives, we may assume that $e$ and $e'$ are hyperimaginary. Consider $a_E\in\bdd^u(Ae)\cap\bdd^u(Ae')$. We may assume $a\ind_{Ae}e'$, whence $ae\ind_A e'$. Let $(a_i:i<\omega)$ be a Morley sequence in $\lstp(a/Ae')$. Then $E(a_i,a_j)$ for all $i,j<\omega$. But $a_i\ind_A a_j$ for $i\not=j$, so $\pi(x,a_j)=\tp(a_i/a_j)$ does not fork over $A$, and neither does $\pi(x,a)$. Note that $\pi(x,y)$ implies $E(x,y)$.

Now suppose $a_E\notin\bdd^u(A)$. We can then find a long sequence
$(a'_i:i<\alpha)$ of $A$-conjugates of $a$ such that $\neg E(a'_i,a'_j)$ for
$i\not=j$. By the Erd\H os-Rado theorem and compactness (see
e.g.\ \cite[Proposition 1.6]{ca11}) there is an infinite $A$-indiscernible
sequence $(a''_i:i<\omega)$ whose $2$-type over $A$ is among the $2$-types of
$(a'_i:i<\alpha)$. In particular $\neg E(a''_i,a''_j)$ for $i\not=j$, and
$(\pi(x,a''_i):i<\omega)$ is $2$-inconsistent. Since $a''_0\models\tp(a/A)$, we
see that $\pi(x,a)$ divides over $A$, a contradiction.\qed

As we have seen in Remark \ref{finchar}, finite character may fail for ultraimaginaries. The next definition singles out the subclass of ultraimaginaries where this does not happen, at least for hyperimaginary sets.
\defn\label{deftame} Let $T$ be simple. An ultraimaginary $e$ is {\em tame} if for all sets $A,B$ of hyperimaginaries we have $e\ind_AB$ if and only if $e\ind_AB_0$ for all finite subsets $B_0\subseteq B$. It is {\em supersimple} if it has a representative of ordinal $SU$-rank.\edefn
\bem\label{quasi-finitary} A supersimple ultraimaginary in a simple theory is
quasi-finitary; in a supersimple theory the converse holds as well.\ebem
\bew Suppose $A$ is a representative for an ultraimaginary $e$ with $SU(A)<\infty$, and let $B$
be a real tuple with $A\in\bdd(B)$. Let $b\in B$ be a finite subtuple with $SU(A/b)$ minimal; it follows that $A\ind_bB$. Hence $A\subseteq\bdd(b)$ and $e$ is bounded over $b$, so $e$ is quasi-finitary. In a supersimple theory the converse is obvious. \qed

We are really interested in the set of tame ultraimaginaries. However, we do not
have a good criterion when an ultraimaginary is tame; moreover, an
ultraimaginary definable over a tame ultraimaginary need not be tame itself. For
instance, the sequence $I$ in Example \ref{bspBY} is tame (since it is real),
but $I_E$ is not. Clearly, an ultraimaginary definable (or even bounded) over a
quasi-finitary/supersimple ultraimaginary is itself quasi-finitary/supersimple.
\lmm\label{tame} A supersimple ultraimaginary is tame. In particular,
quasi-finitary ultraimaginaries in a supersimple theory are tame.\elmm
\bew Let $e$ be a supersimple ultraimaginary, and $a$ a representative with $SU(a)<\infty$. Consider sets $A$ and $B$. There is a finite $b\in B$ with $a\ind_{Ab}B$. So $e\ind_AB$ if and only if $e\ind_Ab$ by Fact \ref{transitive}. Thus $e$ is tame.\qed

In a supersimple theory quasi-finitary ultraimaginaries are the correct ones to consider: Due to elimination of hyperimaginaries all parameters consist of imaginaries of ordinal $SU$-rank; as canonical bases of such imaginaries are finite, we can always reduce to a quasi-finitary situation.

Another kind of tame ultraimaginaries arose in the generalization of the group configuration theorem to simple theories \cite{bytw02,bytw05}.
\defn An invariant equivalence relation $E$ is {\em almost type-definable} if
there is a type-definable symmetric and reflexive relation $R$ finer than $E$
such that any $E$-class can be covered by boundedly many $R$-balls (i.e.\
sets of the form $\{x:xRa\}$ for varying $a$). A class modulo an almost
type-definable equivalence relation is called an {\em almost
hyperimaginary}.\edefn
\tats\cite[page 188]{bytw05} Almost hyperimaginaries are tame. In fact, they satisfy finite character.\etats
We shall now consider how to obtain invariant equivalence
relations, and hence ultraimaginaries.
\satzli\label{geneqrelstab} Let $T$ be stable. For algebraically closed $A$ and an $\emptyset$-invariant equivalence relation $E$ on $\tp(b)$, consider the relation $R(X,Y)$ given by
$$\exists xy\,[Xx\equiv Yy\equiv Ab\land x\ind_XY\land y\ind_YX\land E(x,y)].$$
Then $R$ is an $\emptyset$-invariant equivalence relation on $\tp(A)$.\esatzli
\bew Clearly, $R$ is $\emptyset$-invariant, reflexive and symmetric. So suppose that $R(A,A')$ and $R(A',A'')$ both hold, and let this be witnessed by $b,b'$ and $b^*,b''$. Let $b_1\models\tp(b'/A')=\tp(b^*/A')$ with $b_1\ind_{A'}AA''$. Since $A'$ is algebraically closed, $b'\ind_{A'}A$ and $b^*\ind_{A'}A''$ we have $b_1\equiv_{AA'}b'$ and $b_1\equiv_{A'A''}b^*$ by stationarity. Hence there are $b_0,b_2$ with $bb'\equiv_{AA'}b_0b_1$ and $b^*b''\equiv_{A'A''}b_1b_2$. In particular $E(b_0,b_1)$ and $E(b_1,b_2)$ hold, and so does $E(b_0,b_2)$. Moreover, we may assume $b_0\ind_{AA'b_1}A''$ and $b_2\ind_{A'A''b_1}A'$. Now $b_1\ind_{A'}AA''$ implies $b_0\ind_{AA'}A''$ and $b_2\ind_{A'A''}A$. Then $b_0\ind_AA'$ and $b_2\ind_{A''}A'$ imply $b_0\ind_AA''$ and $b_2\ind_{A''}A$, whence $R(A,A'')$ holds. So $R$ is transitive.\qed

\bem Recall that a reflexive and symmetric binary relation $R(x,y)$ on a partial type $\pi(x)$ is {\em generically transitive} if whenever $x,y,z\models\pi$ and $x\ind_yz$, then $R(x,y)$ and $R(y,z)$ together imply $R(x,z)$. If $T$ is merely simple, the relation $R$ in Proposition \ref{geneqrelstab} is still generically transitive. However, contrary to the type-definable case \cite[Lemma 3.3.1]{wa00}, the two-step iterate of an $\emptyset$-invariant, reflexive, symmetric and generically transitive relation on a Lascar strong type need not be transitive.\ebem
\bsp Consider on the forest of Example \ref{forest} the relation $R(a,b)$ which holds if $3$ divides the distance between $a$ and $b$. This is generically transitive, as for $a'\ind_a a''$ the distance between $a'$ and $a''$ is the sum of the distances between $a'$ and $a$ and between $a$ and $a''$. However, two points of distance $2$ are easily seen to be $R^2$-related, so the transitive closure $E$ of $R$ is just the relation of being in the same connected component. But no two points of distance $1$ are $R^2$-related.\ebsp
Clearly, in the above example the three-step iterate $R^3$ is transitive, as it
is just the relation of being connected. Is there an example of a generically
transitive symmetric and reflexive $\emptyset$-invariant relation $R$ such that
$R^n$ is not transitive for any $n<\omega$~?

The next proposition shows that in a simple theory and under some conditions, if
$R$ is a generically transitive reflexive and symmetric relation, then at least
its transitive closure is not bounded, unless $R$ holds for two independent
elements. We first recall the definitions of $SU_p$-rank and $p$-closure.
\defn\cite[Remark 5.1.19]{wa00} Let $P$ be an $\emptyset$-invariant family of
regular types closed under nonforking extensions. The {\em $SU_P$-rank} is the
smallest function from the collection of all types to the ordinals together with
infinity, such that
$SU(a/A)\ge\alpha+1$ if there is some $B\supseteq A$ and some $c\nind_Ba$ with $\tp(c/B)\in P$.\newline
The {\em $P$-closure} of a set $A$ is given by $\cl_P(A)=\{a:SU_P(a/A)=0\}$.\edefn
Then $SU_P$-rank satisfies the Lascar inequalities \cite[Exercise 5.1.20]{wa00}. Note that unless $P$ contains all non-orthogonality classes of types of $SU$-rank $1$, the $P$-closure has the size of the monster model.
\tats\label{tats1}\cite[Lemma 3.5.3]{wa00} The following are equivalent:\begin{enumerate}
\item $\tp(a/A)$ is foreign to all types $q$ with $SU_P(q)=0$.
\item $a\ind_A\cl_P(A)$.
\item $a\ind_A\dcl(aA)\cap\cl_P(A)$.
\item $\dcl(aA)\cap\cl_P(A)\subseteq\bdd(A)$.\end{enumerate}\etats
Note that $P$-closure is well-behaved with respect to independence:
\tats\label{tats2}\cite[Lemma 3.5.5]{wa00} {\rm and} \cite[Lemma 3]{wa04} Suppose $A\ind_BC$. Then $\cl_P(A)\ind_{\cl_P(B)}\cl_P(C)$.
More precisely, for any $A_0\subseteq\cl_P(A)$ we have
$A_0\ind_{B_0}\cl_P(C)$, where
$B_0=\dcl(A_0B)\cap\cl_P(B)$. In particular,
$\cl_P(AB)\cap\cl_P(BC)=\cl_P(B)$.\etats

\bem\label{general} If $p\in S(\emptyset)$ is regular and $P$ is the family of all non-forking extensions of $p$, we shall write $SU_p$ and $\cl_p$. Another
choice for $P$ is the family of all regular types of $SU$-rank $\omega^\alpha$
used in the proof of Theorem
\ref{EUI}. One can also take $P$ to be the family of types foreign to some
$\emptyset$-invariant collection $\Sigma$ of partial types; in this case
$P$-closure $\cl_P$ is equal to the $\Sigma$-closure
$\cl_\Sigma$ defined in Definition \ref{Sclosure} (see \cite{wa04}).\ebem
\lmm\label{equalrank} Let $P$ be an $\emptyset$-invariant family of regular
types closed under nonforking extensions. Suppose
$SU_P(a/bc)=SU_P(a/b)=n$ is finite. Then $\cl_P(a)\ind_{\cl_P(b)}\cl_P(c)$.\elmm
\bew Let $(a_i:i\le\omega)$ be a Morley sequence in $\lstp(a/bc)$ with
$a=a_\omega$, and put $c_k=(a_i:i<n)$ and $d_k=\cb(a/bc_k)$. We shall show
first that $d_k\in\cl_P(b)$.

Let $(a'_i:i<\omega)$ be a Morley sequence in $\lstp(a/bc_k)$. Since
$$SU_P(d_k/b)\le SU_P(c_k/b)\le nk$$ 
is finite, there is $\ell<\omega$ and
$r\ge0$ with $SU_P(d_k/b,a'_i:i<\ell')=r$ for all $\ell'\ge\ell$. Suppose $r>0$.
Then there is $B\supseteq\{b,a'_i:i<\ell\}$ and $b'$ with $\tp(b'/B)\in P$ such
that $b'\nind_Bd_k$; we
may assume 
$$Bb'\ind_{b,d_k,(a'_i:i<\ell)}(a'_i:i<\omega).$$
As $d_k\in\dcl(a'_i:i<\omega)$, there is minimal $\ell'\ge\ell$ such that
$b'\nind_B(a'_i:i\le\ell')$. So
$$\begin{aligned}SU_P(a'_i:i\le\ell'/b,d_k,
a'_i:i<\ell)&=SU_P(a'_i:i\le\ell'/Bb')\\
&<SU_P(a_i:i\le\ell'/B)\\
&\le SU_P(a'_i:i\le\ell '/b,a'_i:i<\ell)\\
&\le\ell n.\end{aligned}$$
By Lascar symmetry (the second Lascar inequality),
$$r=SU_P(d_k/b,a'_i:i\le\ell')<SU_P(d_k/b,a'_i:i<\ell)=r,$$
a contradiction. Thus $r=0$. By the Lascar inequalites,
$$\begin{aligned}\ell n+SU_P(d_k/b)&=SU_P(a_i:i<\ell/b,d_k)+SU_P(d_k/b)\\
&\le SU_P(d_k,a_i:i<\ell/b)\\
&\le SU_P(d_k/b,a_i:i<\ell)\oplus SU_P(a_i:i<\ell/b)\\
&=0\oplus\ell n,\end{aligned}$$
whence $SU_P(d_k/b)=0$ and $d_k\in\cl_P(b)$.

Since $a\ind_{bd_k}(a_i:i<k)$ Fact \ref{tats2} yields
$a\ind_{\cl_P(b)}(a_i:i<k)$ for all $k$, whence
$a\ind_{\cl_P(b)}(a_i:i<\omega)$. If $d=\cb(a/bc)$, then $d\in\dcl(a_i:i<\omega)$ and $a\ind_{bd}c$. So $\cl_P(a)\ind_{\cl_P(b)}\cl_P(bd)$ and $\cl_P(a)\ind_{\cl_P(bd)}\cl_P(c)$ by Fact \ref{tats2}; the result follows by transitivity.\qed

\satzli\label{geneqrel} Let $T$ be simple. Suppose $R$ is an
$\emptyset$-invariant, reflexive, symmetric and generically transitive relation
on $\lstp(a)$, and $P$ is an $\emptyset$-invariant family of regular types
closed under non-forking extensions such that $SU_P(a)$ is finite. Let $E$ be
the transitive closure of $R$, and suppose $a_E\in\bdd^u(\cl_P(\emptyset))$.
Then there is $a'\ind_{\cl_P(\emptyset)} a$ with $R(a,a')$.\esatzli
\bew Put $c=\bdd(a)\cap\cl_P(\emptyset)$. Then $a\ind_c\cl_P(\emptyset)$ by Fact
\ref{tats1}, whence $a_E\in\bdd^u(c)$ by Lemma \ref{indepintersection}. Let
$a'\equiv_c^{lstp} a$ with $a'\ind_c a$. Then $a_E=a'_E$, so there is $n<\omega$
and a chain $a=a_0,a_1,\ldots,a_n=a'$ such that $R(a_i,a_{i+1})$ holds for all
$i<n$. Put $a'_0=a_0$, and for $0<i<n$ let
$$a'_i\equiv_{a_i}^{lstp}a'_{i-1}\quad\text{with}\quad a'_i\ind_{a_i}a_{i+1}.$$
\beh $\bdd^u(a'_i)\cap\bdd^u(a_{i+1})\subseteq\bdd^u(a_0)$.\ebeh
\bewbeh For $i=0$ this is trivial. For $i>0$, as $a'_i\equiv_{a_i}^{lstp}a'_{i-1}$ and $\bdd^u(a_i)=\dcl^u(\bdd(a_i))$, we get 
$$\bdd^u(a'_i)\cap\bdd^u(a_i)=\bdd^u(a'_{i-1})\cap\bdd^u(a_i).$$
Next, $a'_i\ind_{a_i}a_{i+1}$ implies 
$$\bdd^u(a'_ia_i)\cap\bdd^u(a_ia_{i+1})=\bdd^u(a_i)$$ by Lemma \ref{indepintersection}. Hence inductively
$$\begin{aligned}\bdd^u(a'_i)\cap\bdd^u(a_{i+1})
&\subseteq\bdd^u(a'_i)\cap\bdd^u(a_i)\\
&=\bdd^u(a'_{i-1})\cap\bdd^u(a_i)\\
&\subseteq\bdd^u(a_0).\qed\end{aligned}$$
Now by generic transitivity and induction, $R(a'_i,a_{i+1})$ holds for all $i<n$. In particular $R(a'_{n-1},a_n)$ holds, and by Lemma \ref{indepintersection}
$$\bdd^u(a'_{n-1})\cap\bdd^u(a_n)\subseteq\bdd^u(a_0)\cap\bdd^u(a_n)=\bdd^u(c).$$
Choose $a''$ with $R(a'',a'_{n-1})$  such that $SU_P(a''/a'_{n-1})$ is maximal
possible. We may choose it such that $a''\ind_{a'_{n-1}}a_n$. Then
$$\bdd^u(a'')\cap\bdd^u(a_n)\subseteq\bdd^u(a_n)\cap\bdd^u(a'_{n-1})\subseteq\bdd^u(c)$$ and
$$SU_P(a''/a_n)\ge SU_P(a''/a'_{n-1}a_n)=SU_P(a''/a'_{n-1}).$$
Rename $a''a_n$ as $a_1a_2$, and note that
$\bdd^u(a_1)\cap\bdd(a_2)\subseteq\bdd^u(c)$, $c\subseteq\bdd(a_2)$, and
$SU_P(a_1/a_2)$ is maximal possible among tuples $(x,y)$ with $R(x,y)$.
Moreover,
$$\begin{aligned}
SU_P(a_2/a_1)&=SU_P(a_1a_2)-SU_P(a_1)\\
&=SU_P(a_1a_2)-SU_P(a_2)=SU_P(a_1/a_2 ),\end{aligned}$$
so this is also maximal.

Choose $a_3\ind_{a_2}a_1$ with $a_3\equiv_{a_2}^{lstp}a_1$. By generic transitivity $R(a_1,a_3)$ holds. Moreover,
$$SU_P(a_3/a_1)\ge SU_P(a_3/a_1a_2)=SU_P(a_3/a_2),$$
so equality holds. Similarly, 
$$SU_P(a_1/a_3)=SU_P(a_1/a_2a_3)=SU_P(a_1/a_2).$$
Now $SU_P(a_i/a_j)=SU_P(a_i/a_ja_k)$ for $\{i,j,k\}=\{1,2,3\}$ implies by Lemma
\ref{equalrank} that
$$\cl_P(a_i)\ind_{\cl_P(a_j)}\cl_P(a_k).$$
In particular,
$$\cl_P(a_i)\cap\cl_P(a_k)=\cl_P(a_1)\cap\cl_P(a_2)\cap\cl_P(a_3).$$
Let $b=\cl_P(a_1)\cap\cl_P(a_2)\cap\bdd(a_1a_2)$. Then
$\cl_P(a_1)\cap\cl_P(a_2)=\cl_P(b)$ by \cite[Lemma 3.18]{pw}. Let $F(x,y)$ be
the $\emptyset$-invariant equivalence relation on $\lstp(b)$ given by
$\cl_P(x)=\cl_P(y)$. As $b_F$ is fixed by the $\bdd(a_2)$-automorphism moving
$a_1$ to $a_3$ and $a_1\ind_{a_2}a_3$, we get $b_F\in\bdd^u(a_2)$ by Lemma
\ref{indepintersection}. Similarly, considering an $a_3'\ind_{a_1}a_2$ with
$a'_3\equiv_{a_1}^{lstp}a_2$ we obtain $b_F\in\bdd^u(a_1)$, whence
$$b_F\in\bdd^u(a_1)\cap\bdd^u(a_2)\subseteq\bdd^u(c).$$
So if $b'\models lstp(b/c)$ with $b'\ind_c b$, then
$b'_F=b_F$ and
$$\cl_P(b')=\cl_P(b)=\cl_P(c)=\cl_P(\emptyset).$$
But now 
$$\cb(a_3/\cl_P(a_1)\cl_P(a_2))\subseteq\cl_P(a_1)\cap\cl_P(a_2)=\cl_P(b)=\cl_P(
\emptyset),$$
so $a_3\ind_{\cl_P(\emptyset)}a_2$, as required.\qed

We shall illustrate the use of the proposition in Propositions \ref{internal}
and \ref{bddaugment}, whose proof in the hyperimaginary case uses canonical
bases.
From now on, let $\Sigma$ be an $\emptyset$-invariant family of partial types.
Recall first the definitions of internality, analysability and orthogonality for hyperimaginaries.
\defn Let $\pi$ be a partial type over $A$. Then $\pi$ is\begin{itemize}
\item ({\em almost}) {\em $\Sigma$-internal} if for every realization $a$
of $\pi$ there is $B\ind_Aa$ and a tuple $\bar b$ of realizations of types in $\Sigma$
based on $B$, such that $a\in\dcl(B\bar b)$ (or $a\in\bdd(B\bar b)$, respectively).
\item {\em $\Sigma$-analysable} if for any realization $a$ of $\pi$ there is a sequence
$(a_i:i<\alpha)$ such that $\tp(a_i/A,a_j:j<i)$ is $\Sigma$-internal for all $i<\alpha$, and $a\in\bdd(A,a_i:i<\alpha)$.\end{itemize}
Finally, $p\in S(A)$ is {\em orthogonal} to $q\in S(B)$ if for all $C\supseteq AB$, $a\models p$, and $b\models q$ with $a\ind_A C$ and $b\ind_B C$ we have $a\ind_C b$. The type $p$ is {\em orthogonal to $B$} if it is orthogonal to all types over $B$.
\edefn
Note that in the definition of analysability, we may in addition require $a_i\in\dcl(Aa)$ for all $i<\alpha$. We now generalize these notions to ultraimaginaries.
\defn We shall say that an ultraimaginary $e$ is ({\em almost) $\Sigma$-internal}, or is {\em $\Sigma$-analysable}, if it has a representative which is. Similarly, $e$ is {\em orthogonal} over $A$ to some type $p$ if for all $B\ind_A e$ such that $p$ is over $B$ and for any realization $b\models p|B$ we have $e\ind_ABb$.\edefn
\bem This definition does not imply that we define the notion of an analysis of
an ultraimaginary. Moreover, $e$ orthogonal to $p$ over $A$ need not imply
that $e$ has a representative which is orthogonal to $p$. And unless $e$ is
tame, orthogonality of $e$ over $A$ to $p$ need not imply orthogonality
to $p^{(\omega)}$.\ebem

\defn\label{level} For an ordinal $\alpha$ the {\em $\alpha$-th $\Sigma$-level}
of $a$ over $A$ is defined inductively by $\ell_0^\Sigma(a/A)=\bdd(A)$, and for
$\alpha>0$
$$\ell_\alpha^\Sigma(a/A)=\{b\in\bdd(aA):\tp(b/\bigcup_{\beta<\alpha}\ell_\beta(a/A))\text{ is almost $\Sigma$-internal}\}.$$
We shall write $\ell_\infty^\Sigma(a/A)$ for $\bigcup_\alpha\ell_\alpha^\Sigma(a/A)$, i.e.\ the set of all hyperimaginaries $b\in\bdd(aA)$ such that $\tp(b/A)$ is $\Sigma$-analysable.\edefn
\bem So $a\in\ell_\alpha^\Sigma(a/A)$ is and only if $\tp(a/A)$ is $\Sigma$-analysable in $\alpha$ steps.\ebem
\lmm\label{ortho} If $\tp(a/A)$ is $\Sigma$-analysable in $\alpha$ steps for
some ordinal $\alpha$ or $\alpha=\infty$ and $A\ind b$, put
$c=\ell_\alpha^\Sigma(b)$. Then $Aa\ind_cb$.\elmm
\bew $\cb(Aa/b)$ is definable over a Morley sequence $(A_ia_i:i<\omega)$ in $\lstp(Aa/b)$. Then $(A_i:i<\omega)\ind b$ and $\tp(a_i/A_i)$ is $\Sigma$-analysable in $\alpha$ steps for all $i<\omega$. The union of these analyses level-by-level gives us a $\Sigma$-analysis of $\cb(Aa/b)$ over $(A_i:i<\omega)$ in $\alpha$ steps. As $(A_i:i<\omega)\ind\cb(Aa/b)$, we obtain that $\cb(Aa/b)$ is analysable over $\emptyset$ in $\alpha$ steps, and must be contained in $c$. Thus $Aa\ind_cb$.\qed

\kor\label{foreign} If $c=\ell^\Sigma_\infty(b)$, then $\tp(b/c)$ is foreign to
all $\Sigma$-analys\-able types.\ekor
\bew We apply Lemma \ref{ortho} over $c$ for $\alpha=\infty$, noting that
$\ell_\infty^\Sigma(b/c)=\ell_\infty^\Sigma(b)=c$.\qed

The next proposition is well-known for hyperimaginaries in simple theories,
even without the restriction on $SU_p$-rank: If $\tp(a/A)$ is non-orthogonal to
a regular type $p\in S(\emptyset)$, then there is a $p$-internal
$a_0\in\bdd(aA)\setminus\cl_p(A)$: If $B\ind_Aa$ and $b\models p|AB$
with $b\nind_{AB}a$, just take $a_0=\cb(bB/aA)$. In fact,
if we just require $a_0\in\bdd(aA)\setminus\bdd(A)$, we do not even need $p$ to
be regular \cite[Propositon 3.4.14]{wa00}. However, as for ultraimaginary $a$
the canonical base does not make sense, we have to work harder.
\satzli\label{internal} Let $T$ be simple. Suppose $b_E$ is an ultraimaginary
non-orthogonal to some regular type $p\in S(\emptyset)$, and
$SU_p(\ell_1^p(b))<\omega$. Then there is an almost $p$-internal ultraimaginary
$e\in\bdd^u(b_E)\setminus\bdd^u(\cl_p(\emptyset))$. Moreover,
$e\in\bdd^u(\ell_1^p(b))$.\esatzli
\bew Let $c=\ell_1^p(b)$. Define an $\emptyset$-invariant relation $R$ on $\tp(c)$ by
$$R(c',c'')\quad\Leftrightarrow\quad\exists b'b''\ [b'c'\equiv b''c''\equiv bc\land E(b',b'')].$$
This is reflexive and symmetric; moreover for $c'\ind_{c''}c'''$ with $R(c',c'')$ and $R(c'',c''')$ we can find $b',b'',b^*,b'''$ with
$$b'c'\equiv b''c''\equiv b^*c''\equiv b'''c'''\equiv bc,$$
such that $E(b',b'')$ and $E(b^*,b''')$ hold. Since $c''$ is boundedly closed, $b''\equiv_{c''}^{lstp}b^*$; moreover $b''\ind_{c''}c'$ and $b^*\ind_{c''}c'''$ by Lemma \ref{ortho}. By the Independence Theorem we can assume $b''=b^*$, so $E(b',b''')$ and $R(c',c''')$ hold.
Hence $R$ is generically transitive; let $F$ be its transitive closure. The class $c_F$ is clearly almost $p$-internal. Moreover, if $E(b',b)$ holds there is $c'$ with $b'c'\equiv bc$. Thus $F(c',c)$ holds, so $c_F$ is bounded over $b_E$.

Finally, suppose $c_F\in\bdd^u(\cl_p(\emptyset))$. By Proposition \ref{geneqrel}
there is $c'\ind_{\cl_p(\emptyset)} c$ with $R(c',c)$. Hence there are $b',b^*$
with $b'c'\equiv b^*c\equiv bc$ and $\models E(b',b^*)$. Applying a
$c$-automorphism (and moving $c'$), we may assume $b=b^*$. Let $A\ind b$ be some
parameters and $a$ some realization of $p$ over $A$ with $a\nind_Ab_E$; we may
assume $Aa\ind_bb'$, whence $A\ind bb'$. Moreover $b\ind_c Aa$ by Lemma
\ref{ortho}, whence $b'\ind_cAa$. Thus $b'\ind_{\cl_p(c)}Aa$ by Fact
\ref{tats2}. Now $c'\ind_{\cl_p(\emptyset)}c$ yields
$c'\ind_{\cl_p(\emptyset)}\cl_p(c)$, and hence $c'\ind_{\cl_p(\emptyset)}Aa$.
Then $a\ind_A\cl_p(\emptyset)$ implies $a\ind_Ac'$. Now $b'\ind_{c'}Aa$ by Lemma
\ref{ortho}, whence $b'\ind_A a$. As $b_E=b'_E$ we obtain $a\ind_A b_E$, a
contradiction.\qed

\kor\label{supinternal} Let $e$ be a supersimple ultraimaginary. Suppose $e$ is non-orthogonal to some regular type $p$ over some set $B$. Then there is an almost $p$-internal supersimple $e'\in\bdd^{qfu}(Be)\setminus\bdd^{qfu}(\cl_p(B))$.\ekor
\bew Let $a$ be a representative of $e$ with $SU(a)<\infty$ and put $b=\cb(a/B)$. Then $SU(b)<\infty$, as $b$ is bounded over a finite initial segment of a Morley sequence in $\lstp(a/B)$. Now $e\ind_bB$, so $\tp(e/b)$ is non-orthogonal to $p$. Note that $SU_p(\ell_1^p(a/b)/b)$ is finite by supersimplicity.
By Proposition \ref{internal} applied over $b$ there is an almost $p$-internal
ultraimaginary $e'\in\bdd^u(be)\setminus\bdd^u(\cl_p(b))$; moreover
$e'\in\bdd^u(\ell_1^p(a/b))\subseteq\bdd(ab)$. Thus $e'$ is supersimple, almost
$p$-internal over $b$ and thus over $B$; it is quasi-finitary by Remark
\ref{quasi-finitary}.\qed

\satzli\label{bddaugment} Let $T$ be supersimple. If $AB\ind D$ and $\bdd^{qfu}(A)\cap\bdd^{qfu}(B)=\bdd^{qfu}(\emptyset)$, then $\bdd^{qfu}(AD)\cap\bdd^{qfu}(BD)=\bdd^{qfu}(D)$.\esatzli
\bew We may assume that $A$, $B$ and $D$ are boundedly closed. Consider
$$e\in(\bdd^{qfu}(AD)\cap\bdd^{qfu}(BD))\setminus\bdd^{qfu}(D).$$
Let $p$ be a regular type of least $SU$-rank non-orthogonal to $e$ over $D$.
This exists by transitivity since $e$ is tame. By Corollary \ref{supinternal} we
may assume that $e$ is almost $p$-internal of finite $SU_p$-rank over $D$; let
$a'$ be a representative which is almost $p$-internal over $D$. Put
$a=\cb(a'D/A)$. As $a\ind D$ we obtain that $\tp(a)$ is almost $p$-internal;
note that $SU(a)<\infty$. Since $e\ind_{aD}A$, Lemma \ref{indepintersection}
implies that $e\in\bdd^{qfu}(aD)$. So we may assume that $A=\bdd(a)$ and
$SU_p(A)<\omega$. Moreover, we may assume that $D=\bdd(\cb(aa'/D))$ is the
bounded closure of a finite set.

Let $(A_i:i<\omega)$ be a Morley sequence in $\lstp(A/BD)$ with $A_0=A$, and put $B'=\bdd(A_1A_2)$. Then $B'$ is almost $p$-internal of finite $SU_p$-rank. Since $e\in\bdd^{qfu}(AD)\cap\bdd^{qfu}(BD)$ we have $e\in\bdd^{qfu}(A_iD)$ for all $i<\omega$. Let $e'$ be the set of $B'D$-conjugates of $e$, again a quasi-finitary ultraimaginary. Since any $B'D$-conjugate of $e$ is again in 
$$\begin{aligned}\bdd^{qfu}(A_1D)\cap\bdd^{qfu}(A_2D)&=\bdd^{qfu}(BD)\cap\bdd^{qfu}(A_1D)\\
&=\bdd^{qfu}(BD)\cap\bdd^{qfu}(AD),\end{aligned}$$
we have $e'\in\dcl^{qfu}(B'D)\cap\bdd^{qfu}(AD)$. Moreover, $B'\ind_{BD}A$, whence $B'\ind_BA$ and 
$$\bdd^{qfu}(A)\cap\bdd^{qfu}(B')\subseteq\bdd^{qfu}(A)\cap\bdd^{qfu}(B)=\bdd^{qfu}(\emptyset).$$
Choose $A'\equiv_{AD}^{lstp}B'$ with $A'\ind_{AD}B'$. Then 
$$e'\in\dcl^{qfu}(A'D)\cap\dcl^{qfu}(B'D).$$ 
Now, $D\ind_BA$ implies
$D\ind_BAB'$; as $D\ind B$ we get $D\ind ABB'$. Hence $D\ind_AB'$,
whence $A'\ind_AB'$ and
$$\bdd^{qfu}(A')\cap\bdd^{qfu}(B')\subseteq\bdd^{qfu}(A)\cap\bdd^{qfu}(B')=\bdd^{qfu}(\emptyset).$$
We may assume $e'=(A'D)_E$ for some $\emptyset$-invariant equivalence relation $E$. Define a $\emptyset$-invariant reflexive and symmetric relation $R$ on $\lstp(A')$ by
$$R(X,Y)\Leftrightarrow\exists\,Z\ [XZ\equiv YZ\equiv A'D\land Z\ind XY\land E(XZ,YZ)].$$
By the independence theorem, if $A_1\ind_{A_2}A_3$ such that $R(A_1,A_2)$ and $R(A_2,A_3)$ hold, we have $R(A_1,A_3)$. Hence $R$ is generically transitive; let $E'$ be the transitive closure of $R$. Clearly $A'_{E'}$ is quasi-finitary.

Next, consider $A''\equiv_{B'}A'$ with $A''\ind_{B'}A'$. By the independence theorem there is $D'$ with $A'D\equiv_{B'} A'D'\equiv_{B'}A''D'$ and $D'\ind_{B'}A'A''$. Then $D'\ind B'$, whence $D'\ind A'A''$ and $(A'D')_E=(A''D')_E\in\dcl^{qfu}(B'D')$. Therefore $E'(A',A'')$ holds and $A'_{E'}\in\dcl^{qfu}(B')$. Thus
$$A'_{E'}\in\dcl^{qfu}(A')\cap\dcl^{qfu}(B')\subseteq\bdd^{qfu}(\emptyset).$$
By Proposition \ref{geneqrel} there is $A''\ind_{\cl_p(\emptyset)}A'$ with $R(A',A'')$. Let $D'$ witness $R(A',A'')$. Then $D'\equiv_{A'}D$, so we may assume $D'=D$. Since $\cl_p(D)\ind_{\cl_p(\emptyset)}\cl_p(A'A'')$ and $\cl_p(A')\ind_{\cl_p(\emptyset)}\cl_p(A'')$ we obtain 
$$\cl_p(A')\ind_{\cl_p(\emptyset)}\cl_p(A'')\cl_p(D)$$
and hence $A'\ind_{\cl_p(D)}A''$. But now
$$e'=(A'D)_E=(A''D)_E\in\dcl^{qfu}(A'D)\cap\dcl^{qfu}(A''D)\subseteq\bdd^{qfu}(\cl_p(D))$$
by Lemma \ref{indepintersection}. Since $e\in\bdd^{qfu}(e')$, this contradicts non-orthogonality of $e$ to $p$ over $D$.\qed

\bem Again, the proof of the hyperimaginary analogue of Proposition
\ref{bddaugment} for simple theories uses canonical bases and does not
generalize.\ebem

\section{Elimination of ultraimaginaries}\label{sec4}

One cannot avoid the non-tame ultraimaginaries of Example \ref{bspBY} which do not satisfy finite character and hence cannot be eliminated. Similarly, on a type of rank $\omega$ we cannot eliminate the relation of having mutually finite rank over each other (Example \ref{forest}), since the rank over a class modulo such a relation is not defined. We thus content ourselves with elimination of supersimple ultraimaginaries in a simple theory (and in particular of quasi-finitary ultraimaginaries in a supersimple theory) up to rank of lower order of magnitude. This seems to be optimal, given the examples cited.
\defn Let $e$ be ultraimaginary. We shall say that $SU(a/e)<\omega^\alpha$ if for all representatives $b$ of $e$ we have $SU(a/b)<\omega^\alpha$. Conversely, $SU(e/a)<\omega^\alpha$ if there is a representative $b$ with $SU(b/a)<\omega^\alpha$.\edefn
\bem This does not mean that we define the value of $SU(a/e)$ or of $SU(e/a)$. In fact, one might define
$$SU(e/a)=\min\{SU(b/a):\mbox{ $b$ a representative of $e$}\},$$
but this suggests a precision I am not sure exists. 

On the other hand, Example \ref{forest} shows that
$$SU(a/e)=\sup\{SU(a/b):\mbox{ $b$ a representative of $e$}\}$$
is not a good definition, as the rank of a point $a$ over its connected
component $e=a_E$ would be $\omega$, i.e.\ the same as $SU(a)$.\ebem
\lmm Let $e$ be ultraimaginary. $SU(e/a)<\omega^0$ if and only if $e\in\bdd^{u}(a)$, and $SU(a/e)<\omega^0$ if and only if $a\in\bdd(e)$.\elmm
\bew If $b$ is a representative of $e$ with $SU(b/a)<\omega^0$, then $b\in\bdd(a)$, so $e\in\bdd^{u}(a)$. If $e\in\bdd^{u}(a)$, then $e\in\dcl^{u}(\bdd(a))$, so $b=\bdd(a)$ is a representative of $e$ with $SU(b/a)<\omega^0$.

If $a\notin\bdd(e)$, then there are arbitrarily many $e$-conjugates of $a$. Then for any representative $b$ of $e$ there is some $e$-conjugate $a'$ of $a$ which is not in $\bdd(b)$. Let $b'$ be the image of $b$ under an $e$-automorphism mapping $a'$ to $a$. Then $b'$ is a representative of $e$, and $SU(a/b')\ge\omega^0$. On the other hand, if $a\in\bdd(e)$, then $a\in\bdd(b)$ for any representative $b$ of $e$, whence $SU(a/b)<\omega^0$.\qed

\defn An ultraimaginary $e$ can be {\em $\alpha$-eliminated} if there is a representative $a$ with $SU(a/e)<\omega^\alpha$. A supersimple theory has {\em feeble elimination of ultraimaginaries} if for all ordinals $\alpha$, all quasi-finitary ultraimaginaries of rank $<\omega^{\alpha+1}$ can be $\alpha$-eliminated.\edefn
\bem $0$-elimination is usually called {\em weak} elimination; in the presence of imaginaries this equals full elimination. I do not know what the definition of feeble elimination of ultraimaginaries should be in general for simple theories --- but then their whole theory is much more problematic.\ebem
\satz\label{EUI} If $e$ is ultraimaginary with $SU(e)<\omega^{\alpha+1}$, then $e$ can be $\alpha$-eliminated. A supersimple theory has feeble elimination of ultraimaginaries; a supersimple theory of finite rank has elimination of quasi-finitary ultraimaginaries.\esatz
\bew Let $a$ be a representative of $e$ of minimal rank. Since
$SU(e)<\omega^{\alpha+1}$ we have $SU(a)<\omega^{\alpha+1}$. Suppose
$SU(a/e)\ge\omega^\alpha$, so there is some representative $b$ of $e$ with
$SU(a/b)\ge\omega^\alpha$. Let $P$ be the family of regular types of $SU$-rank
$\omega^\alpha$. Then $SU_P(a)<\omega$ and $SU_P(a/b)=n>0$; we choose $b$ such
that $n$ is maximal. Consider
$a'\equiv_b^{lstp}a$ with $a'\ind_ba$. Since $e\in\dcl^u(b)$ we have
$e\in\dcl^u(a')$. By maximality of $n$,
$$SU_P(a/a')\le n=SU_P(a/b)=SU(a/a'b)\le SU_P(a/a'),$$
so equality holds. By lemma \ref{equalrank} we have
$$a\ind_{\cl_P(a')}\cl_P(b).$$
On the other hand, $a\ind_ba'$ implies by the anlogue of Fact \ref{tats2} that
$$a\ind_{\cl_P(b)}\cl_P(a'),$$
so
$$c=\cb(a/\cl_P(b)\cl_P(a'))\subseteq\cl_P(b)\cap\cl_P(a').$$
Then $a\ind_cb$, so $e\in\bdd^u(c)$ by Lemma \ref{indepintersection}. On the
other hand, $c\in\cl_P(a')\cap\cl_P(b)$ implies
$SU(c/a')<\omega^\alpha$, and $SU(c/b)<\omega^\alpha$. Then $SU(a'/c)\ge
SU(a'/cb)\ge\omega^\alpha$ since $SU(a'/b)\ge\omega^\alpha$. It follows that
$$SU(a)=SU(a')\ge SU(c)+\omega^\alpha.$$
In particular $\bdd(c)$ is a representative for $e$ of lower rank, a contradiction.\qed

\bem Let $p$ be a regular type (or type of weight $1$). Then two realizations $a$ and $b$ of $p$ are independent if and only if $\bdd^{qfu}(a)\cap\bdd^{qfu}(b)=\bdd^{qfu}(\emptyset)$: One direction is Lemma \ref{indepintersection}, the other follows from the observation that dependence is an invariant equivalence relation on realizations of $p$. However, this does not hold for all types: By elimination of quasifinite ultraimaginaries, it is in particular false in non one-based theories of finite rank.\ebem

\section{Decomposition}\label{sec5}

In this section we shall give ultraimaginary proofs of some of Chatzidakis'
results from \cite{zoe} around the weak canonical base property, and suitable
generalisations to the supersimple case. As before, $\Sigma$ will be an
$\emptyset$-invariant family of partial types in a simple theory. 

Recall that $a$ and $b$ are domination-equivalent over $A$, denoted $a\domeq_Ab$, if for any $c$ we have $c\ind_Aa\Leftrightarrow c\ind_Ab$. The following lemma is folklore, but we give a proof for completeness.

\lmm\label{domeq}\begin{enumerate}
\item Suppose $a\domeq b$. If $c\ind a$ and $c\ind b$, then $a\domeq_c b$.
\item Suppose $a\domeq_c b$. If $c\ind ab$ then $a\domeq b$.
\item Suppose $a\domeq_c b$. If $\tp(a)$ and $\tp(b)$ are foreign to $\tp(c)$,
then $a\domeq b$.\end{enumerate}\elmm
\bew\begin{enumerate}
\item Consider any $d$ with $d\nind_ca$. Then $cd\nind a$, whence $cd\nind b$. Now $b\ind c$ implies $b\nind_c d$. The converse follows by symmetry.
\item Consider any $d$ with $d\nind a$. Clearly we may assume $d\ind_{ab}c$, whence $abd\ind c$. Since $a\ind c$ we get $d\nind_c a$, whence $d\nind_c b$ and $cd\nind b$. But $c\ind_d b$, so $d\nind b$; the converse follows by symmetry.
\item Consider any $d$ with $d\nind a$. Since $a\ind c$ we get $d\nind_c a$,
whence $d\nind_c b$ and $cd\nind b$. If $b\ind d$, then $b\ind_dc$ by
foreigness, whence $b\ind cd$, a contradiction. So $b\nind d$; the converse
follows by symmetry.\qed\end{enumerate}

For the following definitions, we require the notion of $\Sigma$-closure alluded
to in Remark \ref{general}.
\defn\label{Sclosure} For an ordinal $\alpha$ we put
$$\cl_\Sigma^\alpha(A)=\{a:\tp(a/A)\mbox{ is $\Sigma$-analysable in
$\alpha$ steps}\}.$$
The {\em $\Sigma$-closure} of $A$ is
$\cl_\Sigma(A)=\cl_\Sigma^\infty(A)=\bigcup_\alpha\cl_\Sigma^\alpha(A)$.\edefn
\bem Note that $\ell_\alpha^\Sigma(a/A)=\cl_\Sigma^\alpha(A)\cap\bdd(aA)$. In
particular, $a\ind_{\ell_\alpha^\Sigma(a/A)}\cl_\Sigma^\alpha(A)$ by
Lemma \ref{ortho}.\ebem
\satzli\label{domeqana} Let $a$ and $b$ be domination-equivalent over
$\cl_\Sigma^\alpha(\emptyset)$, where $a$ is
quasi-finite and $\bdd^{qfu}(a)\cap\bdd^{qfu}(b)=\bdd^{qfu}(\emptyset)$. Then
$ab\in\cl_\Sigma^\alpha(\emptyset)$.\esatzli
\bew Note that by Lemma \ref{ortho} the domination-equivalence of $a$ and
$b$ over $\cl_\Sigma^\alpha(\emptyset)$ means that for any $d$ and
$D=\ell_\alpha^\Sigma(abd)$ we have
$$a\ind_Dd\quad\Leftrightarrow\quad b\ind_Dd.$$
Clearly, domination-equivalence over $\cl_\Sigma^\alpha(\emptyset)$ is an
$\emptyset$-invariant equivalence relation $E$ on $\lstp(a)$. Let
$a'\equiv_{b}^{lstp}a$
with $a'\ind_{b}a$. Then $E(a',a)$ holds. But
$$\bdd^{qfu}(a)\cap\bdd^{qfu}(a')\subseteq\bdd^{qfu}(a)\cap\bdd^{qfu}
(b)=\bdd^{qfu}(\emptyset).$$
Hence $(a)_E=(a')_E\in\bdd^{qfu}(\emptyset)$, and there is $a''\ind a$
with $E(a,a'')$. By Lemma \ref{ortho} we have
$a\ind_{\ell_\alpha^\Sigma(a)}a''\ell_\alpha^\Sigma(aa'')$, whence 
$a\ind_{\ell_\alpha^\Sigma(aa'')}a''$. By domination-equivalence,
$a\ind_{\ell_\alpha^\Sigma(aa'')}a$, that is $a\in\ell_\alpha^\Sigma(aa'')$,
whence $a\in\cl_\Sigma^\alpha(\emptyset)$. Similarly,
$b\in\cl_\Sigma^\alpha(\emptyset)$.\qed

\kor\label{bddeqana} Let $A,B,a,b$ be (hyperimaginary) sets, such that $a$ is
quasi-finite, $\bdd^{qfu}(Aa)\cap\bdd^{qfu}(Bb)=\bdd^{qfu}(\emptyset)$, and $a$
and $b$ are interbounded over $AB$. Suppose $AB$ is $\Sigma$-analysable in
$\alpha$ steps for some ordinal $\alpha$ or $\alpha=\infty$. Then $a$ and $b$
are $\Sigma$-analysable in $\alpha$ steps.\ekor
\bew Since $a$ and $b$ are interbounded over $AB$, they are
domination-equivalent over $\cl_\Sigma^\alpha(\emptyset)$. Now apply 
Proposition \ref{domeqana}.\qed

\bem By Theorem \ref{EUI}, if $SU(Aa)$ or $SU(Bb)$ is finite, then
$\bdd(Aa)\cap\bdd(Bb)=\bdd(\emptyset)$ implies
$\bdd^{qfu}(Aa)\cap\bdd^{qfu}(Bb)=\bdd^{qfu}(\emptyset)$, and we recover
\cite[Lemma 1.15 and Lemma 1.22]{zoe} for $\alpha=\infty$ and $\alpha=1$.\ebem

\tats\label{domeqfor}\cite[Theorem 3.4(3)]{pw} Let $\Sigma'$ be an
$\emptyset$-invariant subfamily of $\Sigma$. Suppose $\tp(a)$ is
$\Sigma$-analysable, but foreign to $\Sigma\setminus\Sigma'$. Then $a$ and
$\ell_1^{\Sigma'}(a)$ are domination-equivalent.\etats

\kor\label{Cb1based} Let $A\subseteq\bdd(\cb(B/A))$ consist of quasi-finite hyperimaginaries, with $\bdd^{qfu}(A)\cap\bdd^{qfu}(B)=\bdd^{qfu}(\emptyset)$. If $A$ is $\Sigma$-analysable and $\Sigma'$ is the subset of one-based partial types in $\Sigma$, then $A$ is analysable in $\Sigma\setminus\Sigma'$.\ekor
\bew Suppose $A$ is not analysable in $\Sigma\setminus\Sigma'$. For every finite tuple $\bar a\in A$ put $c_{\bar a}=\cb(B/\bar a)$, and let $C=\bigcup\{c_{\bar a}:\bar a\in A\}$. Then $A\ind_CB$, as for any $\bar a\in A$ and $C$-indiscernible sequence $(B_i:i<\omega)$ in $\tp(B/C)$ the set $\{\pi(\bar x,B_i):i<\omega\}$ is consistent, where $\pi(\bar x,B)=\tp(\bar a/B)$, since $\pi(\bar x,B)$ does not fork over $c_{\bar a}\subseteq C$. So $A\subseteq\bdd(C)$; as $A$ is not analysable in $\Sigma\setminus\Sigma'$, neither is $C$, and there is $\bar a\in A$ such that $c=c_{\bar a}$ is not analysable in $\Sigma\setminus\Sigma'$. Clearly $c\subseteq\bdd(\bar a)$ is quasi-finite and $c=\cb(B/c)$. Replacing $A$ by $c$ we may thus assume that $A$ is quasi-finite.

Let $A'\subseteq\bdd(A)$ and $B'\subseteq\bdd(B)$ be maximally analysable in $\Sigma\setminus\Sigma'$. So $\tp(A/A')$ and $\tp(B/B')$ are foreign to $\Sigma\setminus\Sigma'$, and $A\not\subseteq A'$. Since $A=\cb(B/A)$ we get $A\nind_{A'}B$; as $A\ind_{A'}B'$ by foreignness to $\Sigma\setminus\Sigma'$, we obtain $A\nind_{A'B'}B$. In particular $B\not\subseteq B'$.

By Fact \ref{domeqfor} the first $\Sigma'$-levels $a=\ell_1^{\Sigma'}(A/A')$ and $b=\ell_1^{\Sigma'}(B/B')$ are non-trivial, one-based, and
$$a\domeq_{A'}A\quad\text{and}\quad b\domeq_{B'}B.$$
Since $\tp(Aa/A')$ is foreign to $\Sigma\setminus\Sigma'$, we have $Aa\ind_{A'}B'$, whence $a\domeq_{A'B'}A$ by Lemma \ref{domeq}(1). Similarly $b\domeq_{A'B'}B$. But $A\nind_{A'B'}B$, and thus $a\nind_{A'B'}b$. Let $a_0=\bdd(A'a)\cap\bdd(A'B'b)$ and $b_0=\bdd(B'b)\cap\bdd(A'B'a)$. By one-basedness of $\tp(a/A')$ and $\tp(b/B')$,
$$a\ind_{A'a_0}B'b\qquad\text{and}\qquad b\ind_{B'b_0}A'a.$$
Hence
$$A'B'a\ind_{A'B'a_0}b_0\qquad\text{and}\qquad A'B'b\ind_{A'B'b_0}a_0.$$
It follows that $a_0$ and $b_0$ are interbounded over $A'B'$.
We can now apply Corollary \ref{bddeqana} to see that $a_0$ is analysable in $\Sigma\setminus\Sigma'$, whence $a_0\in A'$. But then $a\ind_{A'B'}b$, a contradiction.\qed
\bem In a theory of finite $SU$-rank, due to weak elimination of quasi-finitary ultraimaginaries, we obtain that for any $A,B$
$$\tp(\cb(A/B)/\bdd(A)\cap\bdd(B))$$
is analysable in the collection of non one-based types of $SU$-rank $1$.\ebem
\bem Without the quasi-finite hypothesis in Proposition \ref{domeqana},
Corollary \ref{bddeqana} and Corollary \ref{Cb1based}, the conclusions still
hold if we assume that the full ultraimaginary bounded closures intersect
trivially.\ebem
The following Theorem generalizes \cite[Proposition 1.16]{zoe} to supersimple
theories of infinite rank, at the price of demanding that the quasifinite
ultraimaginary bounded closures intersect trivially, rather than just the
bounded closures. The proof is essentially the same, but we have to work with
ultraimaginaries at key steps. Of course, in finite rank this is equivalent, due
to elimination of quasifinite hyperimaginaries; moreover, the families
$\Sigma_i$ in the Theorem are just different orthogonality classes of regular
types of rank $1$.
\defn Two $\emptyset$-invariant families $\Sigma$ and $\Sigma'$ are {\em
perpendicular} if no realization of a type in $\Sigma$ can fork with a
realisation of a type in~$\Sigma'$.\edefn
\bsp If $p$ and $p'$ are two orthogonal types of $SU$-rank $1$ non-orthogonal to
$\emptyset$ (or whose $\emptyset$-conjugates remain orthogonal), then the
families of $\emptyset$-conjugates of $p$ and of $p'$ are perpendicular.\ebsp
\satz\label{decomp} Let $T$ be supersimple. Suppose $A\subseteq\bdd(\cb(B/A))$
and $B\subseteq\bdd(\cb(A/B))$, with
$\bdd^{qfu}(A)\cap\bdd^{qfu}(B)=\bdd^{qfu}(\emptyset)$. Let $(\Sigma_i:i\in I)$
be a family of pairwise perpendicular $\emptyset$-invariant families of partial
types such that $A$ is analysable in $\bigcup_{i\in I}\Sigma_i$. For $i\in I$
let $A_i$ and $B_i$ be the maximal $\Sigma_i$-analysable subset of $\bdd(A)$ and
$\bdd(B)$, respectively. Then $A\subseteq\bdd(A_i:i<\alpha)$ and
$B\subseteq\bdd(B_i:i<\alpha)$; moreover $A_i=\bdd(\cb(B_i/A))$ and
$B_i=\bdd(\cb(A_i/B))$. If $\Sigma_i$ is one-based, then
$A_i=B_i=\bdd(\emptyset)$.\esatz
\bew Since $\cb(A_i/B)$ is $\tp(A_i)$-analysable and hence $\Sigma_i$-analysable, we have $\cb(A_i/B)\subseteq B_i$; similarly $\cb(B_i/A)\subseteq A_i$. As the families in $(\Sigma_i:i\in I)$ are perpendicular, we obtain
$$(A_i:i\in I)\ind_{(B_i:i\in I)}B\quad\text{and}\quad (B_i:i\in I)\ind_{(A_i:i\in I)}A.$$
Suppose $A\subseteq\bdd(A_i:i\in I)$. Then $B=\cb(A/B)\subseteq\bdd(B_i:i\in I)$; moreover
\begin{align*}\bdd(A)&=\bdd(\cb(B/A))=\bdd(\cb(B_i/A):i\in I)\\
&=\bdd(\cb(B_i/A_i):i\in I)\subseteq\bdd(A_i:i\in I)=\bdd(A)\end{align*}
again by perpendicularily. Hence $\bdd(\cb(B_i/A_i))=A_i$, and similarly $\bdd(\cb(A_i/B_i))=B_i$. But if $\Sigma_i$ is one-based, then
$$B_i=\bdd(\cb(A_i/B_i))\subseteq\bdd(A_i)\cap\bdd(B_i)=\bdd(\emptyset)~;$$
similarly $A_i=\bdd(\emptyset)$.

Put $\bar A=\bdd(A_i:i\in I)$ and $\bar B=\bdd(B_i:i\in I)$. It remains to show that $A\subseteq\bar A$. So suppose not. As in the proof of Corollary \ref{Cb1based} put $c_{\bar a}=\cb(B/\bar a)$ for every finite tuple $\bar a\in A$, and let $C=\bigcup\{c_{\bar a}:\bar a\in A\}$. Then again $A\ind_CB$ and $A\subseteq\bdd(C)$; moreover $c_{\bar a}=\cb(B/c_{\bar a})$. Since $A$ is not contained in $\bar A$, neither is $C$. Hence there is $\bar a\in A$ such that $c=c_{\bar a}\notin\bar A$. As the maximal $\Sigma_i$-analysable subset of $\bdd(c)$ is equal to $\bdd(c)\cap A_i$ we may replace $A$ by $c$ and thus assume that $A$ is quasi-finite. Similarly, we may assume that $B$ is quasi-finite.

Since $A=\cb(B/A)\not\subseteq\bar A$, we have $A\nind_{\bar A}B$; as $A\ind_{\bar A}\bar B$ we obtain 
$A\nind_{\bar A\bar B}B$. Let $(b_j:j<\alpha)$ be an analysis of $B$ over $\bar B$ such that for every $j<\alpha$ the type $\tp(b_j/\bar B,b_\ell:\ell<j)$ is $\Sigma_{i_j}$-analysable for some $i_j\in I$. Let $k$ be minimal with $A\nind_{\bar A\bar B}(b_j:j\le k)$. Then $A\ind_{\bar A}\bar B,(b_j:j<k)$ and $\cb(\bar B,(b_j:j\le k)/A)$ is almost $\Sigma_{i_k}$-internal over $\bar A$. Put $A'=\ell_1^{\Sigma_{i_k}}(A/\bar A)$ and $B'=\ell_1^{\Sigma_{i_k}}(B/\bar B)$. Then $A'\not\subseteq\bar A$, and $\cb(A'/B)\subseteq B'$ since $\bar A\ind_{\bar B}B$. Similarly $\cb(B'/A)\subseteq A'$. Moreover $A'\nind_{\bar A\bar B}B$, whence $A'\nind_{\bar A\bar B}B'$. Replacing $A$ by $\cb(B'/A)=\cb(B'/A')$ and $B$ by $\cb(A'/B)=\cb(A'/B')$ we may assume that $\tp(A/\bar A)$ and $\tp(B/\bar B)$ are both almost $\Sigma_k$-internal (where we write $k$ instead of $i_k$ for ease of notation).

\beh $\bdd^{qfu}(AB_k)\cap\bdd^{qfu}(B)=\bdd^{qfu}(B_k)$.\ebeh
\bewbeh Suppose not. As $B$ is analysable in $\bigcup_{i\in I}\Sigma_i$, Corollary \ref{supinternal} yields some $i\in I$ and
$$d\in(\bdd^{qfu}(AB_k)\cap\bdd^{qfu}(B))\setminus\bdd^{qfu}(B_k)$$
such that $d$ is almost $\Sigma_i$-internal over $B_k$; since $\tp(B/B_k)$ is foreign to $\Sigma_k$ we have $i\not=k$. Hence $A\ind_{\bar AB_k}d$, whence $d\in\bdd^{qfu}(\bar AB_k)$ by Lemma \ref{indepintersection}. But $\bar A=\bdd(A_i:i\in I)$ and $d\ind_{A_iB_k}\bar A$ by almost $\Sigma_i$-internality of $d$ over $B_k$, whence $d\in\bdd^{qfu}(A_iB_k)$. If $B_kd\ind A_i$, then $d\ind_{B_k}A_i$ and $d\in\bdd^{qfu}{B_k}$ by Lemma \ref{indepintersection}, contradicting the choice of $d$. Therefore  $B_kd\nind A_i$; by Corollary \ref{supinternal} there is almost $\Sigma_i$-internal $$d'\in\bdd^{qfu}(B_kd)\setminus\bdd^{qfu}(\emptyset).$$
Note that $d'\in\bdd^{qfu}(A_iB_k)\cap\bdd^{qfu}(B)$. But then $d'A_i\ind B_k$, whence $d'\ind_{A_i}B_k$ and $$d'\in\bdd^{qfu}(A_i)\cap\bdd^{qfu}(B)=\bdd^{qfu}(\emptyset),$$
a contradiction.\qed
\beh We may assume $B_k=\bdd(\emptyset)$.\ebeh
\bewbeh Put $A'=\cb(B/AB_k)$. Then $B_k\subset A'=\cb(B/A')$, and $\bdd(A')^{qfu}\cap\bdd(B)^{qfu}=\bdd^{qfu}(B_k)$. If $B'=\cb(A'/B)=\cb(A'/B')$, then $A'\ind_{B'}B$ and $A\ind_{A'}B$ yield $B\ind_{B'}A$ by transitivity, since $B'\subseteq \bdd(B)$. Thus $B\subset\bdd(B')$. We add $B_k$ to the language; note that $B_k\not=\bdd(\emptyset)$ implies $B\nind B_k$, whence $SU(B'/B_k)<SU(B)$. By induction it is thus sufficient to show that $A',B'$ is still a counterexample over $B_k$.

So suppose not, and let $\bdd(A')=\bdd(A'_i:i\in I)$ and $\bdd(B')=\bdd(B_i':i\in I)$ be decompositions, where $A'_i$ and $B'_i$ are maximally $\Sigma_i$-analysable over $B_k$ in $\bdd(A')$ and $\bdd(B')$, respectively. So $B'_k$ is $\Sigma_k$-analysable, whence $B'_k=B_k\subseteq\bar B$ by maximality. Since $B\subseteq\bdd(B')$ is almost $\Sigma_k$-internal over $\bar B$ and $(B'_i:i\not=k)$ is foreign to $\Sigma_k$, we get $B\ind_{\bar B}(B'_i:i\not=k)$, whence $B\subset\bar B$, a contradiction.\qed

By symmetry, we may also assume $A_k=\bdd(\emptyset)$.

Put $B'=\cb(B/A\bar B)$. Then $\bar B\subseteq\bdd(B')$, and since $B$ is almost
$\Sigma_k$-internal over $\bar B$, so is $B'$. If $A'=\cb(B'/A)$, then
$B'\ind_{A'}A$ and $A\ind_{B'}B$ yield $A\ind_{A'}B$, since
$A'\subseteq\bdd(A)$.
Thus $A\subseteq\bdd(A')$. Put $B''=\cb(A/B')=\cb(A/B'')$. Then 
$$B''\subseteq\bdd(B')\subseteq\bdd(A\bar B),$$
and $B''$ is almost $\Sigma_k$-internal over $\bar B$. Moreover, $\bar
A\ind_{\bar B}B'$, whence
$$\bar B=\cb(\bar A/\bar B)=\cb(\bar A/B')\subseteq\cb(A/B')=B''.$$
Finally, $A\ind_{B''}B'$
implies $$A\subseteq\bdd(\cb(B'/A))\subseteq\bdd(\cb(B''/A)).$$
\beh $\bdd^{qfu}(A)\cap\bdd^{qfu}(B'')=\bdd^{qfu}(\emptyset)$.\ebeh
\bewbeh Suppose not. By Corollary \ref{supinternal} there is $i\in I$ and
$$d\in(\bdd^{qfu}(A)\cap\bdd^{qfu}(B''))\setminus\bdd^{qfu}(\emptyset)$$
which is almost $\Sigma_i$-internal; since $A$ is foreign to $\Sigma_k$ we have
$i\not=k$. As $B''$ is almost $\Sigma_k$-internal over $\bar B$ we have
$d\ind_{\bar B}B''$, whence
$$d\in\bdd^{qfu}(A)\cap\bdd^{qfu}(\bar
B)\subseteq\bdd^{qfu}(A)\cap\bdd^{qfu}(B)=\bdd^{qfu}(\emptyset),$$
a contradiction.\qed

Thus $A,B''$ is another counterexample; by induction on $SU(AB/\bar A\bar B)$ we
may assume 
$$SU(A/\bar A\bar B)=SU(AB''/\bar A\bar B)=SU(AB/\bar A\bar B).$$ 
Similarly, there is another counterexample $A'',B''$ with $\bar A\subset
A''\subseteq\bdd(\bar AB'')$, whence 
$$SU(B''/\bar A\bar B)=SU(A''B''/\bar A\bar B)=SU(AB''/\bar A\bar
B)=SU(A/\bar A\bar B).$$
But if $SU(X/Y)=SU(Z/Y)$ with $X\in\bdd(YZ)$, then $X\domeq_YZ$. Thus
$B''\domeq_{\bar A\bar B}A$, so $B''$ and $A$ are domination-equivalent over
$\cl_{\bigcup_{i\not=k}\Sigma_i}(\emptyset)$ by perpendicularity and Lemma
\ref{domeq}(1). So $AB''$ is analysable in
$\bigcup_{i\not=k}\Sigma_i$ by Proposition \ref{domeqana}. But
$\tp(A/\bar A)$ is almost $\Sigma_k$-internal, whence foreign to
$\bigcup_{i\not=k}\Sigma_i$, yielding the final contradiction.\qed

\bem If $C\ind AB$, then $\bdd^{qfu}(A)\cap\bdd^{qfu}(B)=\bdd^{qfu}(\emptyset)$ implies $\bdd^{qfu}(AC)\cap\bdd^{qfu}(BC)=\bdd^{qfu}(C)$ by Lemma \ref{bddaugment}. Hence Theorem \ref{decomp} applies over $C$; this can serve to refine the decomposition.\ebem
\bem In the finite rank context, it is easy to achieve the hypothesis of Theorem \ref{decomp}, as it suffices work over $\bdd(A)\cap\bdd(B)$. In general, however, if
$$\bdd^{qfu}(A)\cap\bdd^{qfu}(B)\supsetneq\bdd^{qfu}(\bdd(A)\cap\bdd(B)),$$
there is no hyperimaginary set $C$ with 
$$\bdd^{qfu}(A)\cap\bdd^{qfu}(B)=\bdd^{qfu}(C),$$
as this equality implies $\bdd(C)=\bdd(A)\cap\bdd(B)$. Thus, 
we cannot work over $\bdd^{qfu}(A)\cap\bdd^{qfu}(B)$, which is not eliminable. If $SU(A/\bdd(A)\cap\bdd(B))<\omega^{\alpha+1}$, feeble elimination nevertheless yields
$$\bdd^{qfu}(A)\cap\bdd^{qfu}(B)\subset\bdd^{qfu}(\cl_\alpha(A)\cap\cl_\alpha(B)),$$
so we can work over $\alpha$-closed sets, as is done in \cite[Theorem 5.4]{pw}.\ebem

\kor Let $T$ be supersimple, and $\Sigma_1$ and $\Sigma_2$ two perpendicular $\emptyset$-invariant families of partial types. Suppose $a$ is quasi-finite, $\tp(a)$ is analysable in $\Sigma_1\cup\Sigma_2$ and $\tp(a/A)$ is $\Sigma_1$-analysable, with $\bdd^{qfu}(a)\cap\bdd^{qfu}(A)=\bdd^{qfu}(\emptyset)$. Then $\tp(a)$ is $\Sigma_1$-analysable.\ekor
\bew Clearly we may assume that $A=\cb(a/A)$. If $a'=\cb(A/a)$, then
$A$ is interbounded with $\cb(a'/A)$. Moreover, as $\tp(a/a')$ is $\Sigma_1$-analysable, $\tp(a')$ is $\Sigma_1$-analysable if and only if $\tp(a)$ is. So we may assume in addition that $a=\cb(A/a)$. 

By Theorem \ref{decomp} we have
$$\bdd(a)=\bdd(\ell_\infty^{\Sigma_1}(a),\ell_\infty^{\Sigma_2}(a)).$$
Hence $\tp(\ell_\infty^{\Sigma_2}(a)/A)$ is $\Sigma_1$-analysable. By perpendicularity, 
$$\ell_\infty^{\Sigma_2}(a)\in\bdd(A)\cap\bdd(a)=\bdd(\emptyset).$$
Hence $a\in\ell_\infty^{\Sigma_1}(a)$ is $\Sigma_1$-analysable.\qed

For $SU(a)$ finite, this specialises to \cite[Proposition 1.20]{zoe}

\end{document}